\documentclass[12pt,reqno]{amsart}

\usepackage{amssymb,amsthm}
\usepackage{amsfonts,cmtiup,comment,stmaryrd}

\usepackage{mathrsfs,cmtiup}

\makeatletter
\def\blfootnote{\xdef\@thefnmark{}\@footnotetext}
\makeatother

\newtheorem{theorem}{Theorem}[section]

\newtheorem{lemma}[theorem]{Lemma}
\newtheorem{proposition}[theorem]{Proposition}
\newtheorem{corollary}[theorem]{Corollary}

\newtheorem{hyp}[theorem]{Hypothesis}
\newtheorem{conj}[theorem]{Conjecture}

\theoremstyle{definition}

\newtheorem{remark}[theorem]{Remark}
\newtheorem*{definition*}{Definition}

\newcommand{\e}{\varepsilon }
\newcommand{\f}{\varphi}

\newcommand{\g}{\gamma}
\newcommand{\F}{{\Bbb F}}

\newcommand{\Z}{{\Bbb Z}}

\newcommand{\bt}{\begin{theorem}}
\newcommand{\et}{\end{theorem}}
\newcommand{\bc}{\begin{corollary}}
\newcommand{\ec}{\end{corollary}}
\newcommand{\bpr}{\begin{proposition}}
\newcommand{\epr}{\end{proposition}}
\newcommand{\be}{\begin{equation}}
\newcommand{\ee}{\end{equation}}
\newcommand{\bp}{\begin{proof}}
\newcommand{\ep}{\end{proof}}
\newcommand{\bconj}{\begin{conj}}
\newcommand{\econj}{\end{conj}}
\newcommand{\bl}{\begin{lemma}}
\newcommand{\el}{\end{lemma}}
\newcommand{\bh}{\begin{hyp}}
\newcommand{\eh}{\end{hyp}}
\newcommand{\br}{\begin{remark}}
\newcommand{\er}{\end{remark}}

\let\leq=\leqslant
\let\geq=\geqslant
\numberwithin{equation}{section}
\newcommand{\ed}{\end{document}}

\begin{document}
\title[Profinite groups with an automorphism]{Profinite groups with an automorphism\\ whose fixed points are right Engel}

\author{C. Acciarri}
\address{Department of Mathematics, University of Brasilia, DF~70910-900, Brazil}
\email{C.Acciarri@mat.unb.br}

\author{E. I. Khukhro}
\address{Charlotte Scott Research Centre for Algebra,\newline\indent University of Lincoln, Lincoln, LN6 7TS, U.K., and\newline\indent  Sobolev  Institute of Mathematics, Novosibirsk, 630090, Russia}
\email{khukhro@yahoo.co.uk}

\author{P. Shumyatsky}
\address{Department of Mathematics, University of Brasilia, DF~70910-900, Brazil}
\email{p.shumyatsky@mat.unb.br}


\begin{abstract}
 An element $g$ of a group $G$ is said to be right Engel if for every $x\in  G$ there is a number $n=n(g,x)$ such that $[g,{}_{n}x]=1$.
We prove that if a profinite group $G$ admits a coprime automorphism $\varphi$ of prime order such that every fixed point of $\varphi$ is a right Engel element, then $G$ is locally nilpotent.
\end{abstract}
\maketitle

\section{Introduction}

Let $G$ be a profinite group, and $\f$ a (continuous) automorphism of $G$ of finite order. We say for short that $\f$ is a \textit{coprime automorphism} of $G$ if its order is coprime to the orders of elements of $G$ (understood as Steinitz numbers), in other words, if $G$ is an inverse limit of finite groups of order coprime to the order of $\f$. Coprime automorphisms of profinite groups have many properties similar to the properties of coprime automorphisms of finite groups. In particular, if $\f$ is a coprime automorphism of $G$, then for any (closed) normal $\f$-invariant subgroup $N$ the fixed points of the induced automorphism (which we denote by the same letter) in $G/N$ are images of the fixed points in $G$, that is, $C_{G/N}(\f )=C_G(\f )N/N$.
Therefore, if $\f$ is a  coprime automorphism  of prime order $q$ such that $C_G(\f)=1$, Thompson's theorem  \cite{tho} implies that $G$ is pronilpotent, and  Higman's theorem~\cite{hig}  implies that $G$ is nilpotent of class bounded in terms of $q$.

In this paper we consider profinite groups admitting a coprime automorphism of prime order all of whose  fixed points are right Engel elements. Recall  that
the $n$-Engel word $[y,{}_{n}x]$ is defined recursively by
$[y,{}_{0}x]=y$ and $[y,{}_{i+1}x]=[[y,{}_{i}x],x]$. An element $g$ of a group $G$ is said to be right Engel if for any $x\in G$ there is an integer $n=n(g,x)$ such that $[g,{}_nx]=1$.  If all elements of a group are right Engel (therefore also left Engel), then the group is called an  Engel group. By a theorem of Wilson and Zelmanov \cite{wi-ze} based on Zelmanov's results \cite{ze92,ze95,ze17} on Engel Lie algebras, an Engel  profinite group is locally nilpotent. Recall that a group is said to be locally nilpotent if every finite subset generates a nilpotent subgroup. In our main result the right Engel condition is imposed on the fixed points of a coprime automorphism  of prime order.

\bt\label{t-main}
Suppose that $\f$ is a coprime automorphism of prime order of a profinite group $G$. If every element of $C_G(\f )$ is a right Engel element of $G$, then $G$ is locally nilpotent.
\et

The proof of Theorem~\ref{t-main} begins with the observation that  a group $G$ satisfying the hypothesis is pronilpotent. Indeed, right Engel elements of a finite group are contained in the hypercentre by the well-known theorem of Baer \cite{ba}. Therefore every finite quotient of $G$ by a $\f$-invariant open normal subgroup is nilpotent by Thompson's theorem~\cite{tho}, since $\f$ acts fixed-point-freely on the quotient by the hypercentre. Assuming in addition that $G$ is finitely generated, it remains to prove that all Sylow $p $-subgroups $S_p $ of $G$ are nilpotent with a uniform upper bound for the nilpotency class. This is achieved in two stages. First  a bound for the nilpotency class of $S_p $ depending on $p $ is obtained for all $p $. Then a bound independent of $p $ is obtained for all sufficiently large primes $p $. At both stages we apply Lie ring methods and the crucial tool is Zelmanov's theorem \cite{ze92,ze95,ze17} on Lie algebras and some of its consequences. Other important ingredients include  criteria for a pro-$p $ group to be $p $-adic analytic in terms of the associated Lie algebra due to Lazard~\cite{laz}, and in terms of bounds for the rank due to Lubotzky and Mann \cite{lu-ma}, and a theorem of Bahturin and Zaicev  \cite{ba-za} on Lie algebras admitting a group of automorphisms whose fixed-point subalgebra is PI.

\section{Preliminaries}

\subsection*{\bf Lie rings and algebras} Products in Lie rings and algebras are called commutators. We use simple commutator notation for left-normed commutators  $[x_1,\dots  , x_k] = [...[x_1, x_2],\dots ,  x_k]$, and the short-hand for Engel commutators $[x,{}_ny]=[x,y,y,\dots, y]$ with $y$ occurring $n$ times.
An element $a$ of a Lie ring or a Lie algebra $L$ is said to be {ad-nilpotent} if there exists a positive integer
$n$ such that $[x,{}_n a] = 0$ for all $x\in L$.
If $n$ is the least integer with this property, then we say that $a$ is ad-nilpotent of index $n$.

The next theorem is a deep result of Zelmanov \cite{ze92,ze95,ze17}.

\bt\label{t-z}
 Let $L$ be a Lie algebra over a field and suppose that $L$ satisfies
a polynomial identity. If $L$ can be generated by a finite set $X$ such that every
commutator in elements of $X$ is ad-nilpotent, then $L$ is nilpotent.
\et

An important criterion for a Lie algebra to satisfy a polynomial identity
is provided by the next theorem, which was proved by Bahturin and Zaicev for
soluble group of automorphisms \cite{ba-za} (and later extended by Linchenko to the general case \cite{lin}). We use the centralizer notation for the fixed point subring $C_L(A)$ of a group of automorphisms $A$ of $L$.

\bt\label{t-bz}
Let $L$ be a Lie algebra over a field $K$. Assume that a finite
group $A$ acts on $L$ by automorphisms in such a manner that $C_L(A)$ satisfies
a polynomial identity. Assume further that the characteristic of $K$ is either $0$
or coprime with the order of $A$. Then $L$ satisfies a polynomial identity.
\et

Both Theorems~\ref{t-z} and \ref{t-bz} admit respective quantitative versions (see for
example \cite{shu00}). For our purposes, we shall need the following proposition for Lie rings proved in \cite{danilo}, which combines both versions. As usual, $\gamma _i(L)$ denotes the $i$-th
term of the lower central series of $L$.

\bpr\label{p-danilo}
Let $L$ be a Lie ring and $A$ a finite group of automorphisms of
$L$ such that $C_L(A)$ satisfies a polynomial identity $f\equiv 0$. Suppose that
$L$ is generated by an $A$-invariant set of $m$ elements such that every commutator in these elements is ad-nilpotent of index at most $n$. Then there exist positive
integers $e$ and $c$ depending only on $|A|$, $f$,  $m$, and $n$ such that $e\gamma _c(L) = 0$.
\epr

We also quote the following useful result proved in \cite[Lemma~5]{kh-sh99} (although it was stated for Lie algebras in  \cite{kh-sh99}, the proof is the same for Lie rings).

\bl\label{l-be}
Let $L$ be a Lie ring, and $M$ a subring of $L$ generated by $m$ elements
such that all commutators in these elements are ad-nilpotent in $L$ of index at
most $n$. If $M$ is nilpotent of class $c$, then for some number $\e =\e (m,n,c)$ bounded in terms of $m$, $n$, $c$ we have $[L,\,\underbrace{M, M,\dots ,M}_\e ]=0$.
\el

\paragraph{\bf Associated Lie rings and algebras}
We now remind the reader of one of the ways of associating a
Lie ring with a group. A series of subgroups of a group $G$
\be\label{e-ser}
G = G_1\geq G_2\geq \cdots
\ee
is called a \textit{filtration} (or an \textit{$N$-series}, or a \textit{strongly central series}) if
\be\label{e-fil}
[G_i,G_j] \leq G_{i+j}\qquad \text{for all}\quad i, j.
 \ee
 For any filtration \eqref{e-ser} we can define an associated Lie ring $L(G)$ with additive group
$$
L(G)=\bigoplus _{i}G_i/G_{i+1},
$$
where the  factors $L_i = G_i/G_{i+1}$ are additively written. The Lie product is defined on
homogeneous elements $xG_{i+1}\in L_i$, $yG_{j+1}\in L_j$ via the group commutators by
$$
[xG_{i+1},\, yG_{j+1}] = [x, y]G_{i+j+1}\in L_{i+j}
$$
and extended to arbitrary elements of $L(G)$ by linearity. Condition~\eqref{e-fil} ensures that this Lie product is well-defined, and group commutator identities imply that $L(G)$ with these operations is a Lie ring. If all factors $G_i/G_{i+1}$ of a filtration \eqref{e-ser} have prime exponent~$p $,
then $L(G)$ can be viewed as a Lie algebra over the field of $p $ elements $\F_p $. If all terms of \eqref{e-ser} are invariant under an automorphism $\f$ of the group $ G$, then $\f$  naturally induces an automorphism of $L(G)$.

We shall normally indicate which filtration is used for constructing an associated Lie ring. One example of a filtration \eqref{e-ser} is given by the lower central series, the terms of which are denoted by $\g _1(G)=G$ and $\g _{i+1}(G)=[\g _i(G),G]$.  It is worth noting that the corresponding associated Lie ring $L(G)$ is generated by the homogeneous component $L_1=G/\g _2(G)$.

 Another example, for a fixed prime number $p $, is the \textit{Zassenhaus $p $-filtration} (also called the \textit{$p $-dimension series}), which is defined by
$$
G_i=\langle g^{p ^k}\mid g\in \g _j(G),\;\, jp ^k\geq i\rangle .
$$
 The factors of this filtration are elementary abelian $p $-groups, so the corresponding associated Lie ring $D_p (G)$ is a Lie algebra over $\F _p $. We denote  by $L_p (G)$ the subalgebra generated by the first factor $G/G_2$. It is well known that the homogeneous components of $D_p (G)$ of degree $s$ coincide with the homogeneous components of  $ L_p (G)$ for all $s$ that are not divisible by $p $. In particular, $L_p (G)$ is nilpotent if and only if $D_p (G)$ is nilpotent. (Sometimes, the notation $L_p (G)$ is used for $D_p (G)$.)

A group $G$ is said to satisfy a \textit{coset identity} if there is a group word $w(x_1,\dots ,x_m)$ and cosets $a_1H,\dots ,a_mH$ of a subgroup $H\leq G$ of finite index such that $w(a_1h,\dots ,a_mh)=1$ for any $h\in H$. Wilson and Zelmanov \cite{wi-ze} proved that if a group $G$ satisfies a coset identity, then the Lie algebra $L_p (G)$ constructed with respect to the Zassenhaus $p $-filtration satisfies a polynomial identity. In fact, the proof of Theorem~1 in \cite{wi-ze} can be slightly modified to become  valid for any filtration $\eqref{e-ser}$ with abelian factors of prime exponent $p $ and the corresponding associated Lie algebra.

\subsection*{\bf Profinite groups} We always consider a profinite group as a topological group. A~subgroup of a topological group will always mean a closed subgroup, all homomorphisms are continuous, and quotients are by closed normal subgroups. This also applies to taking commutator subgroups, normal closures, subgroups generated by subsets, etc. We say that a subgroup is generated by a subset $X$ if it is generated by $X$ as a topological group. Note that if $\f$ is a continuous automorphism of a topological group $G$, then the fixed-point subgroup $C_G(\f )$  is closed. 

Recall that a {pronilpotent} group is a pro-(finite nilpotent) group, that is, an inverse limit of finite nilpotent groups. For a prime $p $, a pro-$p $ group is  an inverse limit of finite $p $-groups. The Frattini subgroup $P'P^p $  of a pro-$p $ group $P$ is the product of the derived subgroup $P'$ and the subgroup generated by all $p $-th powers of elements of $P$. A subset  generates  $P$ (as a topological group) if and only if its image  generates the elementary abelian quotient $P/(P'P^p )$. See, for example, \cite{wil} for these and other properties of profinite groups.

\section{Local nilpotency of Sylow $p $-subgroups}

In this section we prove the local nilpotency of a pro-$p $ group satisfying the hypotheses of the main Theorem~\ref{t-main}. We shall use without special references the fact that fixed points  $C_{G/N}(\f)$ of an automorphism $\f$ of finite coprime order in a quotient by a $\f$-invariant normal open subgroup $N$ are covered by the fixed points in the group: $C_{G/N}(\f)=C_{G}(\f)N/N$.

\bt\label{t-q}
Let $p $ be a prime and suppose that a finitely generated pro-$p $ group $G$ admits  an automorphism $\f$ of prime order $q\ne p $. If every element of $C_G(\f )$ is a right Engel element of $G$, then $G$ is nilpotent.
\et

We begin with constructing a normal subgroup with nilpotent quotient that will  be
the main focus of the proof. Recall that $h(q)$ is a function bounding the nilpotency class of a nilpotent group admitting a fixed-point-free automorphism of prime order $q$ by Higman's theorem~\cite{hig}.

 \bl\label{l-closure}
 There is a finite set $S\subseteq C_G(\f) $ of fixed points of $\f$ such that the quotient $G/H$ by its normal closure $H=\langle S^G\rangle $ is nilpotent of class $h(q)$.
 \el

 \bp
 In the nilpotent quotient $G/\g _{h(q)+2}(G)$ of the finitely generated group $G$ every subgroup is finitely generated. Therefore there is a finite set $S$ of elements of $C_G(\f )$ whose images cover all fixed points of $\f$ in $G/\g _{h(q)+2}(G)$. Let $H=\langle S^G\rangle $ be the normal closure of~$S$. Then the quotient of  $G$ by $H \g _{h(q)+2}(G)$ is nilpotent of class $h(q)$ by Higman's theorem, which means that $\g _{h(q)+1}(G)\leq H  \g _{h(q)+2}(G)$. Since the group $G/H$ is pronilpotent, it follows that $\g _{h(q)+1}(G)\leq H$, as required.
 \ep

 We fix the notation for the subgroup $H=\langle S^G\rangle$ and the finite set  $S\subseteq C_G(\f) $ given by Lemma~\ref{l-closure}. We aim at an application of Zelmanov's Theorem~\ref{t-z} to the associated Lie algebra of $H$, verifying the requisite conditions in a number of steps. The first step is to show that the quotient $G/H'$ is nilpotent, which is achieved by the following lemma.

\bl\label{l-hnilp}
Suppose that $L$ is a finitely generated pro-$p $ group, $M$ is an abelian normal subgroup equal to the normal closure $M=\langle T^L\rangle$ of a finite set $T$ consisting of right Engel elements of $L$, and $L/M$ is nilpotent. Then $L$ is nilpotent.
\el

\bp
We proceed by induction on the nilpotency class of $L/C_L(M)$. The base of induction is the case where $L/C_L(M)$ is abelian, and the corresponding proof follows from the arguments in the step of induction.

Let $T=\{ t_1,\dots ,t_k\}$.
 Let $Z$ be the inverse image of the centre $Z(L/C_L(M))$ of $L/C_L(M)$ (possibly, $Z=L$ in the base of induction). We claim that $Z$ is nilpotent.
 For any fixed $z\in Z$ there are positive integers $n_{i}$ such that $[t_i,{}_{n_{i}}z]=1$. Set $n=\max_{i}n_{i}$; then $[t_i,{}_{n}z]=1$ for all $i$. Moreover, for any $g\in L$ we have $[t_i^g,{}_{n}z]=[t_i,{}_{n}z]^g=1$ since $[z,g]\in C_L(M)$. Since $M=\langle T^L\rangle$ is abelian, this implies  that $[m,{}_{n}z]=1$ for any finite product $m$ of the elements $t_i^g$, $g\in G$. Since these finite products form a dense subset of $M$, we obtain
 \be\label{e-engz}
 [m,{}_{n}z]=1 \qquad \text{for any}\quad  m\in M.
 \ee

 Since $L/M$ is nilpotent and finitely generated, $Z/M$ is nilpotent and finitely generated. Together with \eqref{e-engz} this implies that $Z$ is nilpotent. Indeed, let $Z=\langle M, z_1,\dots ,z_s\rangle$. Any sufficiently long simple commutator in the elements of $M$ and $z_1,\dots ,z_s$ has an initial segment that belongs to $M$ because $Z/M$ is nilpotent. Since $Z/C_Z(M)$ is abelian, the remaining elements (which can all be assumed to be among the $z_i$) can be arbitrarily rearranged without changing the value of the commutator. If the commutator is sufficiently long, one of the $z_i$ will appear sufficiently many times in a row making the commutator trivial by \eqref{e-engz}.

 We now consider $L/Z'$, denoting by the bar the corresponding images of subgroups and elements.   Clearly, $\bar L$, $\bar M$, and $\bar T$ satisfy the hypotheses of the lemma. But now $\bar Z\leq C_{\bar L}(\bar M)$, so the nilpotency class of $\bar L/ C_{\bar L}(\bar M)$ is less than that of $L/C_L(M)$ (unless $Z=L$ when the proof is complete). By the induction hypothesis, $L/Z'$ is nilpotent. Together with the nilpotency of $Z$ proved above, this implies that $L$ is nilpotent by Hall's theorem \cite{hal}.
 \ep

\bl\label{l-fgh}
The subgroup $H$ is generated by finitely many right Engel elements.
\el

\bp
By Lemma~\ref{l-hnilp} applied with $L=G/H'$, $M=H/H'$, and $T=S$, the quotient $G/H'$ is nilpotent. Then $H/H'$ is finitely generated as a subgroup of a finitely generated nilpotent group. The  Frattini quotient $H/(H'H^p )$ is a finite elementary abelian $p $-group. Since $H$ is generated by a set of right Engel elements, conjugates of elements of $C_G(\f )$, we can choose a finite subset of these elements whose images generate $H/(H'H^p )$. Then this finite set also generates the pro-$p $ group $H$.
\ep

 Let $L_p (H)$  be the associated Lie algebra of  $H$ over $\F _p $ constructed with respect to the Zassenhaus $p $-filtration of $H$.

\bpr\label{p-lq}
The Lie algebra $L_p (H)$ is nilpotent.
\epr

\bp
This will follow from Zelmanov's Theorem~\ref{t-z} if we show that $L_p (H)$ satisfies a polynomial identity and is generated by finitely many elements such that all commutators in these elements are ad-nilpotent.

\bl\label{l-pi}
The Lie algebra $L_p (H)$ satisfies a polynomial identity.
\el

\bp
As a profinite Engel group, $ C_H(\f )=H\cap C_G(\f )$ is locally nilpotent by the Wilson--Zelmanov theorem \cite{wi-ze}. It follows that $ C_H(\f )$ satisfies a coset identity on cosets of an open subgroup of $C_H(\f )$. For example, in the group $ C_H(\f )\times C_H(\f )$ the subsets
$$
E_i=\{ (x,y)\in C_H(\f )\times C_H(\f ) \mid [x,{}_iy]=1\}
$$
are closed in the product topology, and
$$
C_H(\f )\times C_H(\f )=\bigcup_{i=1}^{\infty} E_i.
$$
Hence by the Baire category theorem \cite[Theorem~34]{kel}, one of these subsets $E_n$ contains an open subset of  $ C_H(\f )\times C_H(\f )$, which means that there are cosets $x_0K_1$, $y_0K_2$ of open subgroups $K_1,K_2\leq C_H(\f )$ such that $[x,{}_ny]=1$ for all $x\in x_0K_1$ and $y\in y_0K_2$, and therefore for all $x\in x_0(K_1\cap K_2)$ and $y\in y_0(K_1\cap K_2)$. Thus, $ C_H(\f )$ satisfies a coset identity.

The intersections $C_i=C_H(\f )\cap H_i$ with the terms $H_i$ of the Zassenhaus $p $-filtration for $H$ form a filtration of $C_H(\f )$, since obviously, $[C_i, C_j]\leq C_{i+j}$. Let $\hat  L_p (C_H(\f ))$ be the Lie algebra  constructed for $C_H(\f )$ with respect to the filtration $\{ C_i \}$. Since $\f$ is a coprime automorphism, the fixed-point subalgebra $ C_{L_p (H)}(\f )$ is isomorphic to $\hat  L_p (C_H(\f ))$. We apply a version of the Wilson--Zelmanov result \cite[Theorem~1]{wi-ze}, by which a coset identity on $C_H(\f )$ implies that $\hat L_p (C_H(\f ))$ satisfies some polynomial identity. Indeed, the proof of Theorem~1 in \cite{wi-ze} only uses the filtration property $[F_i,F_j]\leq F_{i+j}$ for showing that the homogeneous Lie polynomial constructed from a coset identity on a group $F$ is an  identity of the  Lie algebra constructed with respect to the filtration $\{F_i\}$.

Thus, the fixed-point subalgebra $ C_{L_p (H)}(\f )$ satisfies a polynomial identity. Hence the Lie algebra $L_p (H)$ also satisfies a polynomial identity by the Bahturin--Zaicev Theorem~\ref{t-bz}.
\ep

\bl\label{l-adn}
The Lie algebra $L_p (H)$ is generated by finitely many elements such that all commutators in these elements are ad-nilpotent.
\el

\bp
By Lemma~\ref{l-fgh} the group $H$ is generated by finitely many right Engel elements, say, $h_1,\dots ,h_m$. Their images $\bar h_1,\dots ,\bar h_m$ in the first factor $H/H _2$ of the Zassenhaus $p$-filtration of $H$ generate the Lie algebra $L_p (H)$. Let $\bar c$ be some commutator in these generators $\bar h_i$, and $c$ the same group commutator in the elements $h_i$. For every $j$, since $[h_j,{}_{k_j}c]=1$ for some $k_j=k_j(c)$, we also have $[\bar h_j,{}_{k_j}\bar c]=0$ in $L_p (H)$. We choose a positive integer $s$ such that  $p ^s\geq \max \{k_1,\dots ,k_m\}$. Then $[\bar h_j,{}_{p ^s}\bar c]=0$ for all $j$. In characteristic $p $ this implies that
\be\label{e-qs}
[\varkappa ,{}_{p ^s}\bar c]=0
\ee
for any commutator $\varkappa$ in the $\bar h_i$. This easily follows by induction on the weight of $\varkappa$ from the formula
$$
[[u,v],{}_{p ^s}w]=[[u,{}_{p ^s}w],v]+[u,[v,{}_{p ^s}w]]
$$
that holds in any Lie algebra of characteristic $p $. This formula follows from the Leibnitz formula
$$
[[u,v],{}_{n}w]=\sum_{i=0}^n {n\choose i} [[u,{}_{i}w],\,[v,{}_{n-i}w]]
$$
(where $[a,{}_{0}b]=a$), because the binomial coefficient ${p ^s\choose i}$ is divisible by $p $ unless $i=0$ or $i=p ^s$.

Since any element of  $L_p (H)$ is a linear combination of commutators in the $\bar h_i$, equation \eqref{e-qs} by linearity implies that $\bar c$ is ad-nilpotent of index at most $p ^s$.
\ep

We can now finish the proof of Proposition~\ref{p-lq}.  Lemmas~\ref{l-pi} and \ref{l-adn} show that $L_p (H)$ satisfies the hypotheses of Zelmanov's Theorem~\ref{t-z}, by which $L_p (H)$ is nilpotent. \ep

\bp[Proof of Theorem~\ref{t-q}] By Lemma~\ref{l-closure} the quotient $G/H$ is nilpotent. Being finitely generated, then $G/H$ is a group of finite rank. Here, the rank of a pro-$p $ group is the supremum of the minimum number of (topological) generators over all open subgroups.

The nilpotency of the Lie algebra $L_p (H)$  of the finitely generated pro-$p $ group $H$  established in Proposition~\ref{p-lq} implies that $H$ is a $p $-adic analytic group. This result goes back to Lazard~\cite{laz}; see also \cite[Corollary~D]{sha}. By the Lubotzky--Mann theorem \cite{lu-ma}, being a $p $-adic analytic group is equivalent to being a pro-$p $ group of finite rank.  Thus, both $H$ and  $G/H$ have finite rank, and therefore the whole pro-$p $ group $G$ has finite rank. Hence $G$ is a $p $-adic analytic group and therefore a linear group. By Gruenberg's theorem \cite{grb}, right Engel elements of a linear group  are contained in the hypercentre. Since $H$ is generated by right Engel elements, we obtain that $H$ is contained in the hypercentre of $G$, and since $G/H$ is nilpotent, the whole group $G$ is hypercentral. Being also finitely generated, then $G$ is nilpotent (see \cite[12.2.4]{rob}).
 \ep

\section{Uniform bound for the nilpotency class}

In the main Theorem~\ref{t-main}, we need to prove that if a finitely generated profinite group $G$ admits a coprime automorphism $\f $ of prime order $q$ all of whose fixed points are right Engel in $G$, then $G$ is nilpotent. We already know that $G$ is pronilpotent, and every Sylow $p $-subgroup of $G$ is nilpotent by Theorem~\ref{t-q}. This would imply the nilpotency of $G$ if we had a uniform bound for the nilpotency class of Sylow $p $-subgroups independent of $p $. However, the nilpotency class furnished by the proof of Theorem~\ref{t-q} depends on $p $.

In this section we prove that for large enough  primes $p $ the nilpotency classes of Sylow $p $-subgroups of $G$ are uniformly bounded above in terms of certain parameters of the group $G$. Together with bounds depending on $p $ given by Theorem~\ref{t-q}, this will complete the proof of the nilpotency of~$G$. In the proof, we do not specify the conditions on $p $ beforehand, but proceed with our arguments noting along that our conclusions hold for all large enough primes $p $.

One of the aforementioned parameters is the finite number of generators of $G$, say,~$d$. Clearly, every Sylow $p $-subgroup of $G$ can also be generated by $d$ elements, being a homomorphic image of $G$ by the Cartesian product of all other Sylow subgroups.

\bl\label{l-uni-e}
There are positive integers $n$ and $N_1$ such that  for every $p >N_1$ all  fixed points of $\f$ in the Sylow $p $-subgroup $P$ of $G$ are  right $n$-Engel elements of $P$.
\el

\bp
In the group $C_G(\f )\times G$, the subsets
$$
E_i=\{(x,y)\in C_G(\f )\times G\mid [x,{}_iy]=1\}
$$
are closed in the product  topology. By hypothesis,
$$
\bigcup _iE_i=C_G(\f )\times G.
$$
Hence, by the Baire category theorem \cite[Theorem~34]{kel}, one of these subsets $E_n$ contains an open subset of  $ C_G(\f )\times G$, which means that there are cosets $x_0K$ and  $y_0L$ of open subgroups $K\leq C_G(\f )$ and $L\leq G$ such that $[x,{}_ny]=1$ for all $x\in x_0K$ and $y\in y_0L$. Since
the indices $|C_G(\f ):K|$ and $|G:L|$ are finite, for all large enough primes $p >N_1$ the Sylow $p $-subgroups of $C_G(\f )$ and $G$ are contained in the subgroups $K$ and $L$, respectively. Then for every prime $p >N_1$, in the Sylow $p $-subgroup $P$ the centralizer $C_P(\f )$ consists of right $n$-Engel elements of $P$.
\ep

\bl\label{l-uni-c}
 There are positive integers $c$ and $N_2$ such that for every $p >N_2$ in the Sylow $p $-subgroup $P$ the fixed-point subgroup $C_P(\f )$ is nilpotent of class $c$.
 \el

\bp
By Lemma~\ref{l-uni-e}, for $p >N_1$ in the Sylow $p $-subgroup $P$ the subgroup  $C_P(\f )$ is an $n$-Engel group. By a theorem of Burns and Medvedev \cite{bu-me}, then  $C_P(\f )$ has a normal subgroup $N_p $ of exponent $e(n)$ such that the quotient $C_P(\f )/N_p $ is nilpotent of class $c(n)$, for some numbers $e(n)$ and $c(n)$ depending only on $n$. Clearly, $N_p =1$ for all large enough primes $p >N_2\geq N_1$. Thus, for every prime $p >N_2$ the subgroup $C_P(\f )$ is nilpotent of class $c=c(n)$.
\ep

The following proposition will complete the proof of the main Theorem~\ref{t-main} in view of Lemmas~\ref{l-uni-e} and \ref{l-uni-c}.

\bpr\label{p-unif}
There are functions $N_3(d,q,n,c)$ and $f(d,q,n,c)$ of four positive integer variables $d,q,n,c$ with the following property. Let $p $ be  a prime, and suppose that $P$ is a $d$-generated pro-$p $ group admitting an automorphism $\f$ of prime order $q\ne p $ such that $C_P(\f )$ is nilpotent of class $c$ and consists of right $n$-Engel elements of $P$. If $p >N_3(d,q,n,c)$, then  the group $P$ is nilpotent of class at most  $f(d,q,n,c)$.
\epr

\bp
It is sufficient to obtain a bound for the nilpotency class in terms of $d$, $q$, $n$,~$c$ for every finite quotient $T$ of $P$ by a $\f$-invariant open normal subgroup.
Consider the associated Lie ring $L(T)$ constructed with respect to the filtration consisting of the terms $\g _i(T)$ of the lower central series of $T$:
$$
L(T)=\bigoplus \g _i(T)/\g _{i+1}(T).
$$
As is well known, this Lie ring is nilpotent of exactly the same nilpotency class as $T$ (see, for example, \cite[Theorem~6.9]{khu2}). Therefore it is sufficient to obtain a required bound for the nilpotency class of $L(T)$. We set $L=L(T)$ for brevity. Let $\tilde L=L\otimes _{\Z}\Z [\omega ]$ be the Lie ring obtained by extending the ground ring by a primitive $q$-th root of unity $\omega$. We regard $L$ as $L\otimes 1$ embedded in $\tilde L$. The automorphism of $L$ and of $\tilde L$ induced by $\f$ is denoted by the same letter. Since the order of the automorphism $\f$ is coprime to the orders of elements of the additive group of $\tilde L$, which is a $p $-group, we have the decomposition into analogues of eigenspaces
$$
\tilde L =\bigoplus _{i=0}^{q-1}L_j, \qquad \text{where}\quad L_j=\{x\in \tilde L\mid x^{\f }=\omega ^jx\}.
$$
For clarity we call the additive subgroups $L_j$ \textit{eigenspaces}, and their elements \textit{eigenvectors}. This decomposition can also be viewed as a $(\Z /q\Z)$-grading of $\tilde L$, since
$$
[L_i,L_j]\subseteq L_{i+j\,({\rm mod}\,q)}.
$$
Note that $L_0=C_L(\f )\otimes _{\Z}\Z [\omega ]$.

\bl\label{l-nilp-c}
The fixed-point subring $C_{\tilde L}(\f )$ is nilpotent of class at most $c$.
\el

\bp
Since $\f$ is a coprime automorphism of $T$, we have
$$
C_L(\f )=\bigoplus _i (C_T(\f )\cap \g _i(T))\g _{i+1}(T)/\g _{i+1}(T).
 $$
 Since the fixed-point subgroup $C_T(\f )$ is nilpotent of class $c$, the definition of the Lie products implies that the same is true for $C_L(\f )$ and therefore also for $C_{\tilde L}(\f )=C_L(\f )\otimes _{\Z}\Z [\omega ]$.
\ep

Our main aim is to enable an application to $\tilde L$ of the effective version of Zelmanov's theorem given by Proposition~\ref{p-danilo}. For that we need a $\f$-invariant set of generators of $\tilde L$ such that all commutators in these generators are ad-nilpotent of bounded index.

Let $L_{(k)} =\g _k(T)/\g _{k+1}(T)$ denote the homogeneous component of weight $k$ of the Lie ring $L$, and let $\tilde L _{(k)}=L_{(k)}\otimes _{\Z}\Z [\omega ]$. For clarity we say that elements of  $\tilde L _{(k)}$ or $L _{(k)}$ are \textit{homogeneous}. The component $L_{(1)}$ generates the Lie ring $L$, and $\tilde L_{(1)}$ generates $\tilde L$. If elements $t_1,\dots ,t_d$ generate the group $T$, then their images $\bar t_1,\dots ,\bar t_d$ in $L_{(1)}=T/\g _2(T)$ generate the Lie ring $L$, as well as $\tilde L$ (over the extended ground ring). Writing $\bar t_i=\sum _{j=0}^{q-1}t_{ij}$, where $t_{ij}\in \tilde L_{(1)}\cap L_j$ we obtain a $\f$-invariant set of generators of $\tilde L$
$$
\{\omega ^kt_{ij}\,\mid \,i=1,\dots ,d;\;j=0,\dots ,q-1;\;k=0,\dots ,q-1\}.
$$
 We claim that  for $p >n$ all commutators in these generators are ad-nilpotent of index bounded in terms of $q$, $n$, $c$.

We set for brevity $\tilde L_{(v)k}= \tilde L_{(v)}\cap L_k$ for any weight $v$. A commutator of weight $v$ in the eigenvectors $t_{ij}$ is an eigenvector belonging to $\tilde L_{(v)k}$, where  $k$ is the modulo $q$ sum of the second indices of the $t_{ij}$ involved. We actually prove that  for $p >n$ any homogeneous eigenvector $l_k\in \tilde L_{(v)k}$ is ad-nilpotent of index $s$ bounded in terms of $q$, $n$, $c$. It is clearly sufficient to show that $[x_j,{}_{s}l_k]=0$ for any homogeneous eigenvector $x_j\in \tilde L_{(u)j}$, for any weights $u,v$ and any indices $j,k\in \{0,1,\dots ,q-1\}$. (Here we use indices $j,k$ for elements $x_j,l_k$ only to indicate the eigenspaces they belong to.) First we consider the case where $j=0$.

\bl\label{l-adn0}
If $p >n$, then for any weights $u,v$, for any eigenvector $x_0\in \tilde L_{(u)0}$ and any homogeneous element $l\in \tilde L_{(v)}$ we have $[x_0,{}_{n}l]=0$.
\el

\bp
Since $\f$ is an automorphism of coprime order,  for $x_0\in \tilde L_{(u)0}$ there are elements $y_i\in C_T(\f )\cap \g _u(T)$ such that $x_0=\sum _{i=0}^{q-2}\omega ^i\bar y_i$, where $\bar y_i$ is the image of $y_i$ in $\g _u(T)/\g _{u+1}(T)$ (here the indices of the $y_i$ are used for numbering). For any $\bar h \in L_{(v)}$, there is an element $h\in T\cap \g _v(T)$ such that $\bar h$ is the image of $h$ in $L_{(v)}=\g _v(T)/\g _{v+1}(T)$. Since $[y_i,{}_nh]=1$ in the group $T$ by the hypothesis of Proposition~\ref{p-unif}, we have $[\bar y_i,{}_n\bar h]=0$
in $L$ for every $i$. Hence, by linearity,
\be\label{e-eng}
[x_0,{}_n\bar h]=0
\ee
in $\tilde L$. Note, however, that $\tilde L_{(v)}$ does not consist only of $\Z [\omega ]$-multiples of elements of $L_{(v)}$. Nevertheless, \eqref{e-eng} looks like the $n$-Engel identity, which implies its linearization, which in turn survives extension of the ground ring, and then implies the required property due to the condition $p >n$ making $n!$ an invertible element of the ground ring. However, we cannot simply make a reference to these well-known facts, since this is not exactly an identity, so we reproduce these familiar arguments in our specific situation (jumping over one of the steps).

We substitute  $a_1+\dots +a_n$ for $\bar h$ in \eqref{e-eng} with arbitrary homogeneous elements  $a_i\in  L_{(v)}$  (the indices of the $a_i$ are used for numbering). Thus,
$$
\big[x_0,\,{}_n(a_1+\dots +a_n)\big]=0
$$
for any elements $a_i\in  L_{(v)}$, some of which may also be equal to one another. After expanding  all brackets, we obtain the equation
\be\label{e-mult}
0=\big[x_0 ,\,{}_n(a_1+\dots +a_n)\big]=\sum _{\substack{i_1\geq 0,\dots ,i_n\geq 0 \\ i_1+\dots +i_n=n}}\varkappa _{i_1,\dots ,i_n},
\ee
where $\varkappa _{i_1,\dots ,i_n}$ denotes the sum of all commutators of degree $i_j$ in $a_j$. Replacing $a_1$ with $0$ (only this formal occurrence, keeping  intact all other $a_i$ even if some are equal to $a_1$) shows that
$$
0=\sum _{\substack{i_1= 0,\,i_2\geq 0,\dots ,i_n\geq 0 \\ i_1+\dots +i_n=n}}\varkappa _{i_1,\dots ,i_n}.
$$
Hence we can remove from the right-hand side of \eqref{e-mult} all terms not involving $a_1$  as a formal entry (keeping the other $a_i$ even if some are equal to $a_1$). We  obtain
$$
0=\big[x_0,{}_n(a_1+\dots +a_n)\big]=\sum _{\substack{i_1\geq 1,\,i_2\geq 0,\dots ,i_n\geq 0 \\ i_1+\dots +i_n=n}}\varkappa _{i_1,\dots ,i_n}.
$$
 Then we do the same with $a_2$ for the resulting equation, and so on, consecutively with all the $a_i$. As a result we obtain
$$
0=\sum _{\substack{i_1\geq 1,\dots ,i_n\geq 1 \\ i_1+\dots +i_n=n}}\varkappa _{i_1,\dots ,i_n}=\varkappa _{1,\dots ,1},
$$
that is,
 \be \label{e-mult2}
0=\sum_{\pi\in S_n}\big[x_0,\, a_{\pi (1)},\dots ,a_{\pi (n)}\big],
 \ee
where the right-hand side is the desired linearization.
Every element $l\in \tilde L_{(v)}$ can be written as  a linear combination $l=m_0+\omega m_1+\omega ^2m_2+\dots +\omega ^{q-2}m_{q-2}$, where $m_i\in L_{(v)}$.  Then
\begin{align*}
[x_0, {}_nl]&=\big[x_0 ,\,{}_n(m_0+\omega m_1+\omega ^2m_2+\dots +\omega ^{q-2}m_{q-2})\big]\\ &= \sum _{i=0}^{n(q-2)}\omega ^i\sum _{j_1+2j_2+\dots +(q-2)j_{q-2}=i}\lambda _{j_0,{j_1},\dots ,{j_{q-2}}},
 \end{align*}
where $\lambda _{j_0,{j_1},\dots ,{j_{q-2}}}$ denotes the sum of all commutators in the expansion of the left-hand side with weight $j_s$ in $m_s$. But each of these sums is clearly symmetric and therefore is equal to $0$ as a consequence of \eqref{e-mult2}, where, if an element $a_i$ is required to be repeated $n_i$ times, then the coefficient $n_i!$ appears, which is invertible in the ground ring, since $n_i<p $ and the additive group is a $p $-group. The lemma is proved.
\ep

\bl\label{l-adnk}
If $p >n$, then for any $v$ and $k$, any homogeneous eigenvector $l_k\in \tilde L_{(v)k}$ is ad-nilpotent of index bounded in terms of $q$, $n$, $c$.
\el

\bp
First consider the case $k=0$. Then $l_0=\sum _{i=0}^{q-2}\omega ^i\bar y_i$, where $\bar y_i$ is the image of an element $y_i\in C_T(\f )\cap \g _v(T)$ in $\g _v(T)/\g _{v+1}(T)$ (the indices of the $y_i$ are used for numbering). For each $i$, since $y_i^{-1}$ is a right $n$-Engel element of $T$ by hypothesis, $y_i$ is a left $(n+1)$-Engel element by a result of Heineken~\cite{hei} (see also \cite[12.3.1]{rob}). For any homogeneous element $\bar h\in L_{(u)}$ there is an element $h\in T\cap \g _u(T)$ such that $\bar h$ is the image of $h$ in $\g _u(T)/\g _{u+1}(T)$. Since $[h,{}_{n+1}y_i]=1$ in the group $T$, we have $[\bar h,{}_{n+1}\bar y_i]=0$ in $L$ for every $i$. Hence, by linearity, each $\bar y_i$ is ad-nilpotent in $L$ of index at most $n+1$. Let $M$ be the subring of $L$ generated by $\bar y_0,\bar y_1, \dots ,\bar y_{q-2}$. Since $M\leq C_L(\f )$,
the subring $M$ is nilpotent of class at most $c$ by Lemma~\ref{l-nilp-c}. We can now apply Lemma~\ref{l-be}, by which
$$
[L,\,\underbrace{M, M,\dots ,M}_{\e }]=0
$$
for some $\e =\e (q-1,n+1,c)$ bounded in terms of $q$, $n$, $c$.
This equation remains valid after extension of the ground ring. In particular, $l_0=\sum _{i=0}^{q-2}\omega ^i\bar y_i$ is ad-nilpotent in $\tilde L$ of index at most $\e =\e (q-1,n+1,c)$, as required.

Now suppose that $k\ne 0$. For a homogeneous eigenvector $x_j\in \tilde L_{(u)j}$, the commutator
\be\label{e-subc}
[x_j,{}_{q+n-1}l_k]=\big[x_j,\, \underbrace{l_k,\dots ,l_k}_{s},\,l_k,\dots, l_k\big]
\ee
has an initial segment of length $s+1\leq q$ that is a homogeneous eigenvector $x_0=[x_j,{}_{s}l_k]\in \tilde L_{(w)0}$ (for some weight $w$). Indeed, the congruence $j+sk\equiv 0\,({\rm mod}\, q)$ has a solution $s\in \{0,1,\dots ,q-1\}$ since $k\not\equiv 0\,({\rm mod} \,q)$. There remain at least $n$ further entries of $l_k$ in \eqref{e-subc}, so that we have a subcommutator of the form $[x_0,{}_nl_k]$, which is equal to $0$ by Lemma~\ref{l-adn0}. Thus, by linearity,  $l_k$ is ad-nilpotent of index at most $q+n-1$.
\ep

We now finish the proof of Proposition~\ref{p-unif}. By Lemma~\ref{l-adnk}, for  $p >n$ every commutator in the generators $t_{ij}$ of the Lie ring $\tilde L$ is ad-nilpotent of index bounded in terms of $q$, $n$,~$c$. The same is true for the generators in the  $\f$-invariant set
$$
\{\omega ^kt_{ij}\,\mid \,i=1,\dots ,d;\;j=0,\dots ,q-1;\;k=0,\dots ,q-1\},
$$
which consists of $q^2d$ elements. The fixed-point  subring $C_{\tilde L}(\f )$ is nilpotent of class at most $c$ by Lemma~\ref{l-nilp-c}.
Thus, for $p >n$ the Lie ring $\tilde L$ and its group of automorphisms $\langle \f\rangle$  satisfy the hypotheses of Proposition~\ref{p-danilo}. By this proposition, there exist positive
integers $e$ and $r$ depending only on $d$, $q$, $n$,   $c$ such that $e\gamma _r(\tilde L) = 0$. The additive group of $\tilde L$ is a $p $-group. Therefore,  if $p >e$, then $e$ is invertible in the ground ring, so that we obtain $\gamma _r(L) = 0$. It remains to put $N_3(d,q,n,c)=\max \{n, e\}$ and $f(d,q,n,c)=r-1$.

We thus proved that for $p >N_3(d,q,n,c)$ every finite quotient of $P$ by a $\f$-invariant normal open subgroup is nilpotent of class at most $f(d,q,n,c)$. Therefore $P$ is nilpotent of class at most $f(d,q,n,c)$ if $p >N_3(d,q,n,c)$.
\ep

We finally combine all the results in the proof of the main theorem.

\bp[Proof of Theorem~\ref{t-main}] Recall that $G$ is a profinite group admitting a coprime automorphism $\f $ of prime order $q$ all of whose fixed points are right Engel in $G$; we need to prove that $G$ is locally nilpotent. Any finite set $S\subseteq G$ is contained in the $\f$-invariant finite set $S^{\langle\f\rangle}=\{s^{\f ^k}\mid s\in S,\; k=0,1,\dots ,q-1\}$. Therefore we can assume that the group $G$ is finitely generated, say,  by $d$ elements, and then need to prove that $G$ is nilpotent. As noted in the Introduction, the group $G$ is pronilpotent, so we only need to prove that all Sylow $p $-subgroups of $G$ are  nilpotent of class bounded by  some number independent of~$p $.

Let $n$ and $N_1$ be the numbers given by Lemma~\ref{l-uni-e}, and $c$ and $N_2$ the numbers given by Lemma~\ref{l-uni-c}. Further, let $N_3(d,q,n,c)$ be the number given by  Proposition~\ref{p-unif}. Then for every prime $p >\max\{N_1,N_2,N_3(d,q,n,c)\}$ the Sylow $p $-subgroup of $G$ is nilpotent of class at most $f(d,q,n,c)$ for the function given by  Proposition~\ref{p-unif}. Since every Sylow $p $-subgroup is nilpotent by Theorem~\ref{t-q}, we obtain a required uniform bound for the nilpotency classes of Sylow $p $-subgroups independent of $p $.
\ep

\section*{Acknowledgements}
The first and third authors were supported  by Conselho Nacional de Desenvolvimento Cient\'{\i}fico e Tecnol\'ogico (CNPq), Brazil, and the second   author was supported  by the Russian Science Foundation, project no. 14-21-00065. The second author thanks  CNPq-Brazil and the University of Brasilia for support and hospitality that he enjoyed during his visit to Brasilia.

\end{document}